\documentclass{article}

\title{\LARGE \textbf{Spanning trees in connected graphs with few branch and end vertices}}
\author{Zhora Nikoghosyan\footnote{G.G. Nicoghossian (up to 1997)}  }

\begin{document}

\maketitle

\begin{abstract}

A vertex of degree one in a tree  is called an end vertex and a vertex of degree at least three is called a branch vertex. For a graph $G$, let $\sigma_2$ be the  minimum degree sum of two nonadjacent  vertices in $G$. We consider tree problems arising in the context of optical and centralized terminal networks: finding a spanning tree of G (i) with the minimum number of end vertices, (ii) with the minimum number of branch vertices and  (iii) with the minimum degree sum of the branch vertices, motivated by network design problems where junctions are significantly more expensive than simple end- or through-nodes, and are thus to be avoided. We consider: $(\ast)$ connected graphs on $n$ vertices such that $\sigma_2\ge n-k+1$ for some positive integer $k$. In 1976, it was proved (by the author) that every graph satisfying  $(\ast)$ has a spanning tree with at most $k$ end vertices. In this paper we first show that every graph satisfying $(\ast)$ has a spanning tree with at most $k+1$ branch and end vertices altogether. The next result states that every graph satisfying $(\ast)$ has a spanning tree with at most $(k-1)/2$ branch vertices. The third result states that every graph satisfying $(\ast)$ has a spanning tree with at most $\frac{3}{2}(k-1)$ degree sum of branch vertices. All results are sharp.  \\


\end{abstract}

\section{Introduction}

Throughout this article we consider only finite undirected graphs without loops or multiple edges. The set of vertices of a graph $G$ is denoted by $V(G)$ and the set of edges by $E(G)$.  A good reference for any undefined terms is $\cite{[1]}$.

For a graph $G$, we use $n$ and $\alpha$ to denote the order (the number of vertices) and the independence number of $G$, respectively.  If $\alpha\ge k$ for some integer $k$, let $\sigma_k$ be the minimum degree sum of an independent set of $k$ vertices; otherwise we let $\sigma_k=+ \infty$. For a subset $S\subseteq V(G)$, we denote by $G[S]$ the subgraph of $G$ induced by $S$. We use $d_G(v)$ to denote the number of neighbors of a vertex $v$  in $G$, called the degree of $v$ in $G$.

If $Q$ is a path or a cycle in a graph $G$, then the order of $Q$, denoted by $|Q|$, is $|V(Q)|$. The graph $G$ is hamiltonian if $G$ contains a Hamilton cycle, i.e. a cycle containing every vertex of $G$.

We write a path $Q$ with a given orientation by $\overrightarrow{Q}$. For $x,y\in V(Q)$, we denote by $x\overrightarrow{Q}y$ the subpath of $Q$ in the chosen direction from $x$ to $y$.  We use $x^+$ to denote the successor, and $x^-$ the predecessor, of a vertex $x\in V(Q)$. For $X\subseteq V(Q)$, we define $X^+=\{x^+:x\in X\}$ and $X^-=\{x^-:x\in X\}$.

A vertex of degree one is called an end-vertex, and an end-vertex of a tree is usually called a leaf. The set of end-vertices of $G$ is denoted by $End(G)$. A branch vertex of a tree is a vertex of degree at least three. The set of branch vertices of a tree $T$ will be denoted by $B(T)$.  For a positive integer $k$, a tree $T$ is said to be a $k$-ended tree if $|End(T)|\le k$. A Hamilton path is a spanning 2-ended tree. A Hamilton cycle can be interpreted as a spanning 1-ended tree.

We begin with two famous results on Hamilton paths. \\

\noindent\textbf{Theorem A} \cite{[8]}. Every graph with $\sigma_2\ge n-1$ has a Hamilton path.\\

\noindent\textbf{Theorem B} \cite{[3]}. Every $s$-connected $(s\ge1)$ graph with $\alpha\le s+1$ has a Hamilton path.\\

There are several problems on spanning trees which are generalizations of the Hamilton path problem motivated from optimization aspects with various applications. In this paper we consider tree problems arising in the context of optical networks: (i) finding a spanning tree of G with the minimum number of end vertices, (ii) finding a spanning tree with the minimum number of branch vertices and  (iii) finding a spanning tree of G such that the degree sum of the branch vertices is minimized, motivated by network design problems where junctions are significantly more expensive than simple end- or through-nodes, and are thus to be avoided.

The problem of finding a spanning tree with bounded number of leaves appear when designing centralized terminal networks.

In 1971, Las Vergnas \cite{[6]} gave a degree condition that guarantees that any forest in $G$ of limited size and with a limited number of leaves can be extended to a spanning tree of $G$ with a limited number of leaves in an appropriate sense. This result implies as a corollary a degree sum condition for the existence of a tree with at most $k$ leaves including Theorem A as a special case for $k=1$.   \\

\noindent\textbf{Theorem C} \cite{[2]}, \cite{[6]}, \cite{[7]}. Let $G$ be a connected graph with $\sigma_2\ge n-k+1$ for some positive integer $k$. Then $G$ has a spanning $k$-ended tree.\\

However, Theorem C was first openly formulated and proved in 1976 by the author \cite{[7]}. Later, it was reproved in 1998 by  Broersma and Tuinstra \cite{[2]}.

Win \cite{[9]} obtained a generalization of Theorem B.\\

\noindent\textbf{Theorem D} \cite{[9]}. Let $G$ be a $s$-connected graph with $\alpha\le s+k-1$ for some integer $k\ge2$. Then  $G$ has a spanning $k$-ended tree.\\

One of the interest in the existence of spanning trees with bounded number of branch vertices arises in the realm of multicasting in optical networks.

Gargano, Hammar, Hell, Stacho and Vaccaro \cite{[5]} proved the following.\\

\noindent\textbf{Theorem E} \cite{[5]}. Every connected graph with $\sigma_3\ge n-1$ has a spanning tree with at most one branch vertex.\\

Flandrin et al. \cite{[4]} posed the following conjecture.\\

\noindent\textbf{Conjecture A} \cite{[4]}. If $G$ is a connected graph with $\sigma_{k+3}\ge n-k$ for some positive integer $k$, then $G$ has a spanning tree with at most $k$ branch vertices.\\

In this paper we present a sharp Ore-type condition for the existence of spanning trees in connected graphs with bounded total number of branch and end vertices improving Theorem C by incorporating the number of branch vertices as a parameter. \\

\noindent\textbf{Theorem 1}. Let $G$ be a connected graph of order $n$. If $\sigma_2\ge n-k+1$ for some positive integer $k$, then $G$ has a spanning tree $T$ with at most $k-|B(T)|+1$ end vertices.\\

Let $G$ be the complete bipartite graph $K_{\delta,\delta+k-1}$ of order $n=2\delta+k-1$ and minimum degree $\delta$, where $k\ge3$. Clearly, $\sigma_2(G)=2\delta=n-k+1$. By Theorem 1, $G$ has a spanning tree $T$ with $|End(T)|\le k-b+1$. Observing that $T$ is not $(k-1)$-ended, that is $|End(T)|\ge k$, we have $b\le 1$. On the other hand, we have $b\ge 1$, since $|End(T)|\ge k\ge 3$, which implies $b=1$. This means that $T$ is not $(k-b)$-ended and consequently, Theorem 1 is sharp for each $k\ge3$.\\

The next result follows from Theorem 1 providing a sharp Ore-type condition for the existence of spanning trees in connected graphs with few branch vertices.\\

\noindent\textbf{Theorem 2}. Let $G$ be a connected graph of order $n$. If $\sigma_2\ge n-k+1$ for some positive integer $k$, then $G$ has a spanning tree with at most $(k-1)/2$ branch vertices.\\

The third result provides an Ore-type condition for the existence of spanning trees in connected graphs with bounded degree sum of the branch vertices.\\

\noindent\textbf{Theorem 3}. Let $G$ be a connected graph of order $n$. If $\sigma_2\ge n-k+1$ for some positive integer $k$, then $G$ has a spanning tree with at most $\frac{3}{2}(k-1)$ degree sum of the branch vertices.\\

Let $G$ be a graph (tree) obtained from the path $v_0v_1...v_bv_{b+1}$ by adding new vertices $u_1,...,u_b$ and the edges $u_iv_i$ $(i=1,...,b)$. Clearly, $n=2b+2$ and $\sigma_2=2=n-(2b+1)+1$. Since $|B(G)|=b$, the bound $(k-1)/2$ in Theorem 2 is sharp. Further, since $\sum_{i=1}^bd(v_i)=\frac{3}{2}(k-1)$, the bound $\frac{3}{2}(k-1)$ in Theorem 3 is sharp as well.

\section{Proof of Theorem 1}

\noindent\textbf{Proof of Theorem 1}. Let $G$ be a connected graph with $\sigma_2\ge n-k+1$ and let $T$ be a spanning tree in $G$. Assume that\\

$(a1)$ $T$ is chosen so that $|End(T)|$ is as small as possible.\\

Put $End(T)=\{\xi_1,...,\xi_f\}$. Let $\overrightarrow{P_2}=\xi_1\overrightarrow{P_2}\xi_2$ be the unique path in $T$ with end vertices $\xi_1$ and $\xi_2$. Further, assume that \\

$(a2)$ $T$ is chosen so that $P_2$ is as long as possible, subject to $(a1)$.\\

Put $|B(T)|=b$. If $f=2$ then $P_2$ is a 2-ended spanning tree (Hamilton path) in $G$ with $|B(P_2)|=b=0$, implying that $f=2\le k+1=k-b+1$.

Now let $f\ge3$, that is $b\ge1$. \\

\textbf{Claim 1}. If $P$ is a Hamilton path in $G[V(P_2)]$ with end vertices $x,y$, then $N(x)\cup N(y)\subseteq V(P_2)$.

\textbf{Proof}. Assume the contrary and assume w.l.o.g. that $N(x)\not\subseteq V(P_2)$. Put $T^\prime=T-E(P_2)+E(P)$. Clearly, $T^\prime$ is an $f$-ended spanning tree in $G$ and  $xv\in E(G)$ for some $v\in V(G-P)$. Let $C$ be the unique cycle in $T^\prime+xv$ and let $vv^\prime$ be the unique edge on $C$ with $v^\prime\not= x$. Then $T^\prime+xv-vv^\prime$ is an $f$-ended spanning tree in $G$, contradicting $(a2)$.          \ \ \   $\triangle$ \\

By Claim 1,  $N(\xi_1)\cup N(\xi_2)\subseteq V(P_2)$. If $N(\xi_1)\cap N^+(\xi_2)\not=\emptyset$ then clearly, $G[V(P_2)]$ has a Hamilton cycle. Since $b\ge 1$, $G[V(P_2)]$ has a Hamilton path with end vertex $x$ such that $N(x)\not\subseteq V(P_2)$, contradicting Claim 1. Hence, $N(\xi_1)\cap N^+(\xi_2)=\emptyset$. Observing also that $\xi_1\not\in N(\xi_1)\cup N^+(\xi_2)$ and $N^+(\xi_2)\subseteq V(P_2)$, we get
$$
|P_2|\ge|N(\xi_1)|+|N^+(\xi_2)|+|\{\xi_1\}|
$$
$$
= d(\xi_1)+d(\xi_2)+1\ge \sigma_2+1.          \eqno{(1)}
$$
For each $i\in\{3,...,f\}$, let $\overrightarrow{P_i}=\xi_i\overrightarrow{P_i}z_i$  be the unique path in $T$ between $\xi_i$ and the nearest vertex $z_i$ of $P_2$. Clearly, $z_i\in B(T)$ \ $(i=3,...,f)$.\\

\textbf{Case 1}. $|P_i|=2$ $(i=3,...,f)$.

It follows that $B(T)\subseteq V(P_2)$. If $b=1$ then by (1), $|P_2|\ge \sigma_2+b$ and therefore,
$$
f=|\{\xi_3,...,\xi_f\}|+2=n-|P_2|+2
$$
$$
\le n-\sigma_2-b+2\le k-b+1.
$$

Let $b\ge 2$ and let $x_1,...,x_b$ be the elements of $B(T)$, occurring on $\overrightarrow{P_2}$ in a consecutive order. Assume w.l.o.g. that $x_1=z_3$. Further, assume that \\

$(a3)$ $T$ is chosen so that $d_T(x_1)$ is as small as possible, subject to $(a1)$ and $(a2)$. \\

If $\xi_3v_1\in E(G)$ for some $v_1\in V(x_1^+\overrightarrow{P_2}\xi_2)$, then $T+\xi_3v_1-\xi_3x_1$ is an $f$-ended tree, contradicting (a3). Hence, we can assume that $N(\xi_3)\subseteq V(\xi_1\overrightarrow{P_2}x_1)$, that is
$$
(N(\xi_3)-z_3)\cap B(T)=\emptyset.          \eqno{(2)}
$$
Next, if $N^-(\xi_1)\cap (N(\xi_3)-z_3)$ has an element $v_2$, then
$$
v_2\overleftarrow{P_2}\xi_1v_2^+\overrightarrow{P_2}\xi_2
$$
is a Hamilton path in $G[V(P_2)]$ with end vertex $v_2$ such that $N(v_2)\not\subseteq V(P_2)$, contradicting Claim 1. Hence,
$$
N^-(\xi_1)\cap (N(\xi_3)-z_3)=\emptyset.             \eqno{(3)}
$$
Finally, if $N^-(\xi_1)\cap B(T)\not=\emptyset$, that is $\xi_1z_i^+\in E(G)$ for some $i\in \{3,...,f\}$, then
$$
z_i\overleftarrow{P_2}\xi_1z_i^+\overrightarrow{P_2}\xi_2
$$
is a Hamilton path in $G[V(P_2)]$ with end vertex $z_i$ such that $N(z_i)\not\subseteq V(P_2)$, again contradicting Claim 1. Hence,
$$
N^-(\xi_1)\cap B(T)=\emptyset.              \eqno{(4)}
$$
Using (2), (3), (4) and observing that $\xi_2\not\in N^-(\xi_1)\cup (N(\xi_3)-z_3)\cup B(T)$, we get
$$
|V(P_2)|\ge |N^-(\xi_1)|+|N(\xi_3)-z_3|+|B(T)|+|\{\xi_2\}|
$$
$$
\ge d(\xi_1)+d(\xi_3)+b\ge \sigma_2+b,
$$
implying that
$$
f=|\{\xi_3,...,\xi_f\}|+2=n-|V(P_2)|+2
$$
$$
\le n-\sigma_2-b+2\le k-b+1.
$$

\textbf{Case 2}. $|P_i|\ge3$ for some $i\in\{3,...,f\}$, say $i=3$.

\textbf{Case 2.1}. $N^{-}(\xi_1)\cap N^{+}(\xi_2)\not=\emptyset$.

It follows that $\xi_1w^+,\xi_2w^-\in E(G)$ for some $w\in N^{-}(\xi_1)\cap N^{+}(\xi_2)$. If $z_3=w$ then
$$
w\overleftarrow{P_2}\xi_1w^+\overrightarrow{P_2}\xi_2
$$
is a Hamilton path in $G[V(P_2)]$ with end vertex $w$ such that $N(w)\not\subseteq V(P_2)$, contradicting Claim 1. Hence $z_3\not=w$. Assume w.l.o.g. that $z_3\in V(\xi_1\overrightarrow{P_2}w^-)$. Put
$$
T^\prime=T+\xi_1w^++\xi_2w^--z_3z_3^--ww^-.
$$
Clearly, $T^\prime$ is a spanning $f$-ended tree in $G$ and
$$
\xi_3\overrightarrow{P_3}z_3\overrightarrow{P_2}w^-\xi_2\overleftarrow{P_2}w^+\xi_1\overrightarrow{P_2}z_3^-
$$
is a path in $T^\prime$ longer than $P_2$, contradicting (a2).\\

\textbf{Case 2.2}. $N^{-}(\xi_1)\cap N^{+}(\xi_2)=\emptyset$.

Put
$$
B_1=V(P_2)\cap B(T), \  B_2=B(T)-B_1.
$$
Using Claim 1, it is easy to see that
$$
N^{-}(\xi_1)\cap B_1= N^{+}(\xi_2)\cap B_1=\emptyset.
$$
Observing also that $N^{-}(\xi_1)\cup N^{+}(\xi_2)\subseteq V(P_2)$, we get
$$
|P_2|\ge |N^{-}(\xi_1)|+|N^{+}(\xi_2)|+|B_1|
$$
$$
=d(\xi_1)+d(\xi_2)+|B_1|\ge \sigma_2+|B_1|\ge n-k+1+|B_1|.
$$
Then
$$
n\ge |P_2|+|B_2|+|\{\xi_3,...,\xi_f\}|
$$
$$
\ge n-k+1+|B_1|+|B_2|+f-2=n-k+b+f-1,
$$
implying that $f\le k-b+1$.   \ \qquad     \rule{7pt}{7pt}\\

\noindent\textbf{Proof of Theorem 2}. By Theorem 1, $G$ has a spanning tree $T$ with $|End(T)|\le k-b+1$, where $b=|B(T)|$. On the other hand, it is not hard to see that $|End(T)|\ge b+2$, implying that $b\le (k-1)/2$.   \ \qquad    \rule{7pt}{7pt}\\

\noindent\textbf{Proof of Theorem 3}. By Theorem 1 and Theorem 2, $G$ has a spanning tree $T$ with $f=|End(T)|\le k-b+1$ and $b\le (k-1)/2$, where $b=|B(T)|$. Let $d_1,d_2,...,d_b$ be the degrees of branch vertices of $T$. Observing that
$$
f=\sum_{i=1}^b(d_i-2)+2,
$$
we get 
$$
\sum_{i=1}^bd_i\le k+b-1\le \frac{3}{2}(k-1).
$$

\noindent Institute for Informatics and Automation Problems\\ National Academy of Sciences\\
P. Sevak 1, Yerevan 0014, Armenia\\
E-mail: zhora@ipia.sci.am

\end{document}